\documentclass[a4paper,11pt]{article}
\usepackage{amscd, amsmath, amssymb, amsfonts, epic}
\usepackage{graphicx}
\usepackage{a4wide}
\usepackage{xypic}
\usepackage{paralist}
\usepackage{booktabs}
\usepackage{multirow}
\usepackage{epsfig}
\usepackage{graphics}

\title{The cone conjecture for abelian varieties}
\author{Artie Prendergast-Smith}

\date{}

\def\Z{\text{\bf Z}}
\def\Q{\text{\bf Q}}
\def\R{\text{\bf R}}
\def\C{\text{\bf C}}
\def\H{\text{\bf H}}

\def\arrow{\rightarrow}

\def\iso{\cong}

\def\Pic{\text{Pic}} 
\def\Eff{\text{Eff}}
\def\Aut{\text{Aut}}
\def\End{\text{End}}
\def\Hom{\text{Hom}}

\def\Tr{\text{Tr}}
\def\Amp{A}
\def\PsAut{\text{PsAut}}

\newcommand{\Nef}[1]{\overline{A(#1)}}
\newcommand{\Mov}[1]{\overline{M(#1)}}

\newtheorem{theorem}{Theorem}[section]
\newtheorem{lemma}[theorem]{Lemma}
\newtheorem{corollary}[theorem]{Corollary} 
\newtheorem{proposition}[theorem]{Proposition}

\newtheorem{conjecture}[theorem]{Conjecture}

\begin{document}

\maketitle

The purpose of this paper is to write down a complete proof of the
Morrison--Kawamata cone conjecture for abelian varieties. The
conjecture predicts, roughly speaking, that for a large class of
varieties (including all smooth varieties with numerically trivial
canonical bundle) the automorphism group acts on the nef cone with
rational polyhedral fundamental domain. (See Section \ref{sect-conj}
for a precise statement.) The conjecture has been proved in dimension
2 by Sterk--Looijenga, Namikawa, Kawamata, and Totaro
\cite{Sterk1985,Namikawa1985,Kawamata1997,Totaro2009}, but in higher
dimensions little is known in general.

Abelian varieties provide one setting in which the conjecture is
tractable, because in this case the nef cone and the automorphism
group can both be view as living inside a larger object, namely the
real endomorphism algebra. In this paper we combine this fact with
known results for arithmetic group actions on convex cones to produce
a proof of the conjecture for abelian varieties. It should be stressed
that our proof is a more or less straightforward deduction from
well-known results, and it is written out in full here to serve as a
reference.

Here is the main result.

\begin{theorem} \label{theorem-maintheorem}
Let $X$ be an abelian variety and $\Nef{X}^e$ its effective nef
cone. Then there is a rational polyhedral fundamental domain for the
action of the automorphism group $\Aut(X)$ on $\Nef{X}^e$.
\end{theorem}
The conclusion of the theorem was already known in some cases. It was
proved for abelian surfaces by Kawamata \cite{Kawamata1997}, adapting
the proof of Sterk--Looijenga for $K3$ surfaces. In the same paper,
Kawamata also proved the conjecture for all self-products of an
elliptic curve without complex multiplication. Finally, Bauer
\cite{Bauer1998} showed that for an abelian variety, the nef cone is
rational polyhedral if and only if the variety is isogenous to a
product of mutually non-isogenous abelian varieties of Picard number
1, which in particular implies the theorem for abelian varieties of
this special type.

Thanks to Eugene Eisenstein and Lars Halvard Halle for their comments,
and to Burt Totaro for suggesting the proof presented here.

\section{The cone conjecture} \label{sect-conj}
We work throughout the paper over an arbitrary algebraically closed
field.

Morrison \cite{Morrison1993} gave the original statement of the cone
conjecture for Calabi--Yau threefolds, motivated by considerations
from mirror symmetry. The statement was generalised by Kawamata
\cite{Kawamata1997} to families of varieties with numerically trivial
canonical bundle, and from there to so-called klt Calabi--Yau pairs
\cite{Totaro2009}. As mentioned in the introduction, the conjecture
has been verified in dimension 2, but in general it remains wide
open. See \cite[Section 1]{Totaro2009} for history and a summary of
the current status.

Here we state the conjecture in a rather simple form applicable to
abelian varieties. The symbol $\equiv$ denotes numerical equivalence
of divisors, and for a projective variety $N^1(X)$ denotes the real
vector space $(\text{Div}(X) / \equiv) \otimes \R$ where
$\text{Div}(X)$ is the free abelian group spanned by Cartier
divisors. The cones $\Nef{X}$ and $\Mov{X}$ are the closed cones in
$N^1(X)$ spanned by the classes of nef or movable divisors. The cone
$\Eff(X)$ is the cone spanned by the classes of all effective
divisors, and $\Nef{X}^e$ and $\Mov{X}^e$ denote the intersections
$\Nef{X} \cap \Eff(X)$ and $\Mov{X} \cap \Eff(X)$. Finally, a {\it
  pseudo-automorphism} of $X$ is a birational map $X \dashrightarrow
X$ which is an isomorphism outside a subset of codimension 2. Note
that a pseudo-automorphism maps a movable or effective divisor to
another movable or effective divisor, therefore preserves the cone
$\Mov{X}^e$.

\begin{conjecture}[Morrison--Kawamata]
Let $X$ be a smooth projective variety with $K_X \equiv 0$. Then:

\noindent (1) There exists a rational polyhedral cone $\Pi$ which
is a fundamental domain for the action of $\Aut(X)$ on
$\Nef{X}^e$ in the sense that

\quad (a) $\Nef{X}^e = \Aut(X) \cdot \Pi$,

\quad (b) For all $g \in \Aut(X)$ we have $(\text{Int} \ \Pi) \cap g^*
(\text{Int} \ \Pi) = \emptyset$ for $g^* \neq 1$ in $GL(N^1(X))$.

\noindent (2) There exists a rational polyhedral cone $\Pi'$ which
is a fundamental domain (in the sense above) for the action of the
pseudo-automorphism group $\PsAut(X)$ on $\Mov{X}^e$.
\end{conjecture}
Part (1) of the conjecture would imply in particular that any variety
with numerically trivial canonical bundle has finitely many
contractions up to automorphisms. This is because any
contraction of a projective variety is determined by a semi-ample line
bundle, and two semi-ample line bundles give the same contraction if
they belong to the interior of the same face in the nef effective
cone.

For abelian varieties part (2) of the conjecture is implied by part
(1), because the effective nef cone $\Nef{X}^e$ and effective movable
cone $\Mov{X}^e$ are the same for abelian varieties. Indeed, on an
abelian variety any effective divisor is semi-ample by the 'theorem of
the square' \cite[\S 6, Corollary 4]{Mumford1970}, so
$\Nef{X}^e=\Mov{X}^e=\Eff(X)$. As a consequence Theorem
\ref{theorem-maintheorem} implies the full cone conjecture for abelian
varieties.

The equality of cones in the last paragraph can be strengthened to give
the following result; a reference is
\cite[Proposition 1.1]{Bauer1998}.
\begin{proposition} \label{prop-bauer}
Let $D$ be a Cartier divisor on an abelian variety $X$ and let $[D]$
denote the class of $D$ in $N^1(X)$. Then $[D] \in \Eff(X)$ if and
only if $[D] \in \Nef{X}$.
\end{proposition}
This implies that $\Nef{X}^e$ is equal to the {\it rational hull}
$\Amp(X)_+$ of the ample cone $\Amp(X)$, defined as the convex hull of
the rational points in the closure $\Nef{X}$. This will be important
later, because from the point of view of reduction theory for
arithmetic groups $\Amp(X)_+$ is a more natural object to consider
than $\Nef{X}^e$.

Finally it should be emphasised that for abelian varieties, it makes
no difference to the conjecture whether $\Aut$ denotes the group of
automorphisms in the category $Var/k$ of varieties over the ground
field $k$ or in the category $GrVar/k$ of group varieties over
$k$. This is because any $(Var/k)$-automorphism of an abelian variety
can be composed with a translation automorphism to give a
$(GrVar/k)$-automorphism, and translations act trivially on
$N^1(X)$. For the rest of the paper we will use $\Aut$ to denote the
group of $(GrVar/k)$-automorphisms.

{\bf Examples.} The following examples illustrate how the nef
effective cone and automorphism group of an abelian variety can vary,
compatibly with the cone conjecture. We will consider two abelian
surfaces of Picard number 2. By standard results on quadratic forms
and the Hodge Index Theorem, for any surface $X$ we can choose a
rational basis $\{v_1,v_2\}$ for $N^1(X)$ in which the intersection
form has matrix $\operatorname{diag}(a,-b)$, with $a$ and $b$
positive. When $X$ is an abelian surface, Proposition \ref{prop-bauer}
says that a Cartier divisor $D$ on $X$ is ample if and only if $D^2>0$
and $D \cdot H>0$ for some fixed ample divisor $H$. So if we choose a
basis $\{v_1,v_2\}$ for $N^1(X)$ as above, the ample cone of $X$ is
described as
\begin{align*}
\Amp(X) = \left\{ x_1v_1+x_2v_2 \in N^1(X) \,| \, ax_1^2-bx_2^2 > 0,
\, x_1>0 \right\}. 
\end{align*}
The two extremal rays of $\Nef{X}$ are spanned by the vectors $v_1 \pm
(\sqrt{a/b}) v_2$, so $\Nef{X}^e=\Amp(X)_+$ is a rational polyhedral
cone if and only if $a/b$ is a square in $\Q$.

Before giving our examples we mention the following useful
fact. To verify the cone conjecture for a variety $X$, it suffices to
find a rational polyhedral cone $\Pi \subset \Nef{X}^e$ whose
translates by $\Aut(X)$ cover the whole nef effective cone. This fact
will follow from Theorem \ref{theorem-looijenga} in the next section,
but we mention it now to clarify our examples.

For the first example we take $X$ to be a product $E_1 \times E_2$ of
non-isogenous elliptic curves. Then $\Aut(X) = \Aut(E_1) \times
\Aut(E_2)$, and $\Aut(E_i)$ is a cyclic group of order 2, 4, or
6. Therefore in this case $\Aut(X)$ is a finite group acting on
$\Nef{X}^e$. On the other hand, by taking suitable rational linear
combinations of the divisors $E_1 \times \{0\}$ and $\{0\} \times
E_2$, the intersection form on $N^1(X)$ can be transformed to have
matrix $\operatorname{diag}(1,-1)$. By our description of the extremal
rays a few paragraphs ago, $\Nef{X}^e$ is therefore a rational
polyhedral cone, so by the fact in the previous paragraph the cone
conjecture is true for $X$, taking $\Pi=\Nef{X}^e$.

\begin{figure}[h]
\centerline{
\includegraphics[scale=1]{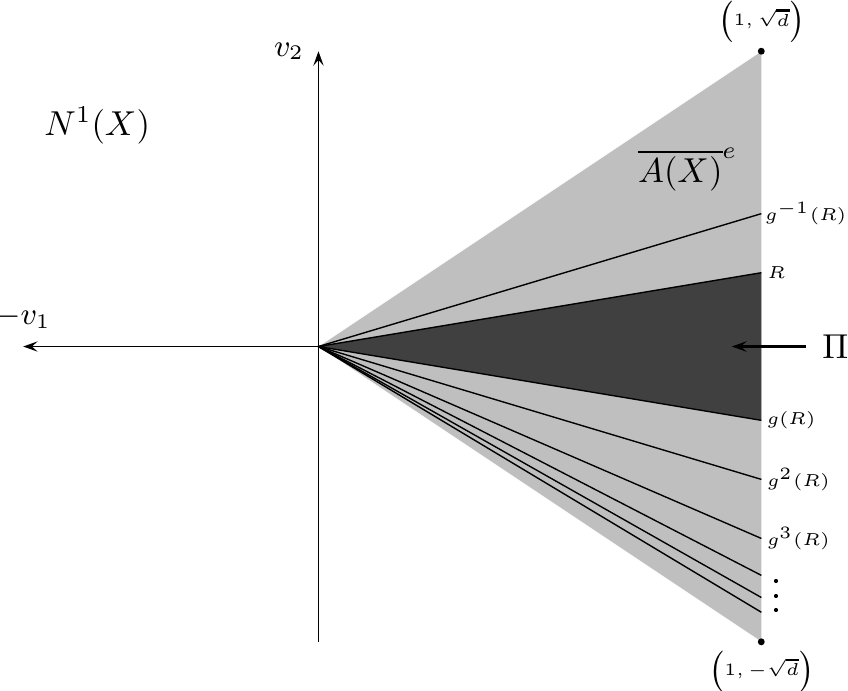}
}
\caption{Nef cone of an abelian surface with real multiplication}
\label{figure-realmult}
\end{figure}

For the second example we take $X$ to be an abelian surface with {\it
  real multiplication}, by which we mean that the endomorphism algebra
$\End^0(X) := \End(X) \otimes \Q$ is isomorphic to the number field
$\Q(\sqrt d)$ for some square-free integer $d>0$. (The simplest
examples of such surfaces are Jacobians of certain genus 2 curves;
explicit models can be found in \cite{Wilson2000}.) By Dirichlet's
unit theorem, the automorphism group $\Aut(X)$ has rank equal to
$r_1+r_2-1$, where $r_1$ is the number of embeddings $\Q(\sqrt d)
\hookrightarrow \R$ and $r_2$ is the number of conjugate pairs of
embeddings $\Q(\sqrt d) \hookrightarrow \C$ whose image is not
contained in $\R$. One checks easily that $r_1=2$ and $r_2=0$, so
$\Aut(X)$ has rank 1. What about the cone $\Nef{X}^e$? In this case
the matrix of the intersection form diagonalises to
$\operatorname{diag}(1,-d)$, so the boundary rays are irrational. To
find the required rational polyhedral fundamental domain, we proceed
as follows. Choose an arbitrary rational ray $R \subset \Nef{X}^e$ and
an element $g$ of infinite order in $\Aut(X)$. A little thought shows
$\left\{ g^i(R) \, | \, i=1,2,3 \ldots \right\}$ is a sequence of
rational rays, which either converges to one extremal ray, or
decomposes into two subsequences, one converging to each extremal
ray. Composing $g$ with a torsion element of $\Aut(X)$, we can assume
that the first case occurs. Then the cone $\Pi$ spanned by the rays $R$
and $g(R)$ is a rational polyhedral cone whose translates by $\Aut(X)$
cover the whole nef effective cone (Figure \ref{figure-realmult}).
Again by the fact above, this proves the cone conjecture for $X$.

\section{Homogeneous self-dual cones and reduction theory} \label{sect-reductionthy}

In this section we give the results we need about reduction theory for
arthimetic group actions on homogeneous self-dual convex cones. It
should be mentioned that this theory has a rich history that we touch
on here only briefly; see \cite{AMRT1975} for historical discussion
and references.

From now on $V$ will always denote a finite-dimensional real vector
space. By a {\it cone} in $V$ we always mean a convex cone $C \subset
V$ which is {\it non-degenerate} (meaning that its closure $\bar{C}$
contains no nonzero subspaces of $V$). The {\it dual cone} $C^* \subset
V^*$ is defined to be the interior of the cone $\overline{C^*}$ consisting
of linear forms on $V$ which are nonnegative on $C$.

Now suppose $C$ is an open cone in a vector space $V$. We define the
{\it automorphism group} $G(C)$ of $C$ to be the subgroup of $GL(V)$
consisting of linear transformations which preserve $C$. The cone $C$
is said to be {\it homogeneous} if $G(C)$ acts transitively on
$C$. Suppose further that $V$ carries an inner product, giving an
identification of $V$ with $V^*$. We say $C$ is {\it self-dual} if
this identification takes $C$ to its dual $C^*$. (This condition
depends on the choice of inner product, but the dependence will not
matter in what follows.) The basic theorem about the automorphism
group of a homogeneous self-dual cone is the following, due to Vinberg
\cite{Vinberg1965}.
\begin{theorem}[Vinberg] \label{theorem-vinberg}
Let $C \subset V$ be a homogeneous self-dual convex cone. Then the
automorphism group $G(C)$ is the group of real points of a reductive
algebraic group $\mathcal{G}(C)$.
\end{theorem}

It is an amazing fact that homogeneous self-dual convex cones can be
completely classified into a small number of cases. More precisely,
define the {\it direct sum} of cones $C_1$ and $C_2$ in vector spaces
$V_1$ and $V_2$ to be the cone $C_1 \oplus C_2 := \{v_1+v_2 \in V_1
\oplus V_2 | v_i \in C_i \}$ and call a cone {\it indecomposable} if
it is not the direct sum of 2 nontrivial cones. The classification
theorem, due to Koecher and Vinberg, is the following \cite[I, \S 1,
  Proposition 2]{Vinberg1963}, \cite[p. 71]{Vinberg1965}.

\begin{theorem}[Koecher---Vinberg] \label{theorem-homogcones}
Any convex cone $C$ can be written as a direct sum $\oplus_i C_i$ of
indecomposable cones. The product $\prod G(C_i)$ is a
finite-index subgroup of $G(C)$. The cones $C_i$ are homogeneous and
self-dual if and only if $C$ is. Any indecomposable homogeneous
self-dual cone is isomorphic to one of the following:
\begin{enumerate}
\item the cone $P_r(\R)$ of positive-definite matrices in the space
  $\mathcal{H}_r(\R)$ of  $r \times r$ real symmetric matrices;
\item the cone $P_r(\C)$ of positive-definite matrices in the space
  $\mathcal{H}_r(\C)$ of  $r \times r$ complex Hermitian matrices;
\item the cone $P_r(\H)$ of positive-definite matrices in the space
  $\mathcal{H}_r(\H)$ of $r \times r$  quaternionic Hermitian matrices;
\item the spherical cone $\left\{ (x_0,\ldots,x_n) \in \R^{n+1} \, |
  \, x_0 > \sqrt{ x_1^2 + \cdots x_n^2} \right\}$; 
\item the 27-dimensional cone of positive-definite $3 \times 3$
  octonionic Hermitian matrices. 
\end{enumerate}
The inner product for which the cone is self-dual is $\left< x,y
\right> = \operatorname{Tr}(xy^*)$ in all cases except 4, and the
usual inner product  on $\R^{n+1}$ in case 4.
\end{theorem}
The proof uses the surprising correspondence between self-dual
homogeneous convex cones and formally real Jordan algebras; see
\cite[Chapter II, \S 2]{AMRT1975} for details.

Vinberg \cite{Vinberg1965} computed the automorphism groups of all the
cones in the list of Theorem \ref{theorem-homogcones}. Here we state
the part of his result which will be relevant to abelian varieties.

\begin{theorem}[Vinberg] \label{theorem-vinberg-auts}
Let $C$ be one of the cones $P_r(k)$ in the previous theorem: that is,
the cone of positive-definite matrices in the vector space
$\mathcal{H}_r(k)$ of $r \times r$ symmetric or Hermitian matrices
over $k$, where $k=\R$, $\C$, or $\H$. Then the identity component
$G(C)^0$ of the automorphism group of $C$ consists of all $\R$-linear
transformations of $\mathcal{H}_r(k)$ of the form $D \mapsto M^* D M$
for some $M \in GL(r,k)$.
\end{theorem}

Now we come to reduction theory. For a $\Q$-algebraic group
$\mathcal{G}$, a subgroup $\Gamma \subset \mathcal{G}(\Q)$ is called
an {\it arithmetic subgroup} of $\mathcal{G}$ if $\Gamma$ and the
group $\mathcal{G}(\Z)$ of integer points of $\mathcal{G}$ are
commensurable (meaning their intersection inside $\mathcal{G}(\Q)$ has
finite index in each). The basic problem of reduction theory for
convex cones is the following: given a convex cone $C$ and an
arithmetic subgroup $\Gamma$ of the automorphism group
$\mathcal{G}(C)$, can we find a rational polyhedral fundamental domain
for the action of $\Gamma$ on $C$? The first results of this kind go
back to Hermite and Minkowski, who found rational polyhedral
fundamental domains for the adjoint action of $SL(r,\Z)$ on the cone
$P_r(\R)$ of positive-definite real symmetric matrices. That is, there
is a finite set of integral linear inequalities such that any
quadratic form can be reduced by an integral change of basis to a form
whose coefficients satisfy those inequalities. This explains the name
`reduction theory'.

The main result of the theory we will use is the following, due to Ash
\cite[Chapter II]{AMRT1975}.

\begin{theorem}[Ash] \label{theorem-ash}
Let $C$ be a self-dual homogeneous convex cone in a real vector space
$V$ with $\Q$-structure. Let $G(C)$ be the automorphism group of $C$
and $\mathcal{G}(C)$ the associated reductive algebraic group (as in
Theorem \ref{theorem-vinberg}). Assume $\mathcal{G}(C)^0$ is defined
over $\Q$. Then for any arithmetic subgroup $\Gamma$ of $\mathcal{G}$,
there exists a rational polyhedral cone $\Pi \subset C_+$ such that
$(\Gamma \cdot \Pi) \cap C = C$.
\end{theorem}
Here as before $C_+$ is the rational hull of $C$, meaning the convex
hull of the rational points in $\overline{C}$.

Applied in the case of abelian varieties, this theorem will provide us
with a rational polyhedral cone whose translates by automorphisms
cover the ample cone. The cone conjecture asks for more, however: we
need to find a precise fundamental domain for the action of
automorphisms on the effective nef cone. This is supplied by the
following result of Looijenga \cite[Proposition 4.1, Application
  4.15]{Looijenga2009}.

\begin{theorem}[Looijenga] \label{theorem-looijenga}
Let $C$ be a convex cone in a vector space $V$ with $\Q$-structure,
and $\Gamma$ a subgroup of $GL(V)$ which preserves $C$ and some
rational lattice in $V$. If there exists a rational polyhedral cone
$\Pi \subset C_+$ such that $(\Gamma \cdot \Pi) \cap C = C$, then in
fact $\Gamma \cdot \Pi = C_+$, and in fact there exists a rational
polyhedral fundamental domain for the action of $\Gamma$ on $C_+$.
\end{theorem}
We observed in Section \ref{sect-conj} that if $X$ is an abelian
variety then $\Amp(X)_+ = \Nef{X}^e$, so together Theorems
\ref{theorem-ash} and \ref{theorem-looijenga} will imply the cone
conjecture for abelian varieties if we can show that for any abelian
variety, the ample cone is a self-dual homogeneous cone and the
automorphism group of the variety acts as an arithmetic subgroup of
the automorphism group of the cone. This is what we do in the next two
sections.

\section{The endomorphism algebra of an abelian variety}

In this section we describe the endomorphism algebra of an abelian
variety as a product of certain matrix algebras. In particular this
gives a description of the automorphism group as an arithmetic group,
connecting the cone conjecture for abelian varieties to the reduction
theory of the previous section. For a full exposition of the structure
theory of the endomorphism algebra, read Chapter 4 of Mumford's
beautiful book \cite{Mumford1970}.

The first result we need is a form of the Poincar\'e complete
reducibility theorem \cite[\S 19, Theorem 1, Corollary
  1]{Mumford1970}.
\begin{theorem} \label{thm-complete-reducubility}
Let $X$ be an abelian variety. Then $X$ is isogenous to a product
$X_1^{n_1} \times \cdots \times X_k^{n_k}$ where the $X_i$ are simple
abelian varieties, not isogeneous for $i \neq j$. The isogeny type of
the $X_i$ and the natural numbers $n_i$ are uniquely determined by
$X$.
\end{theorem}
Here a {\it simple} abelian variety is one that has no proper abelian
subvarieties. An {\it isogeny} $X \arrow Y$ of abelian varieties is a
surjective homomorphism with finite kernel, and $X$ and $Y$ are {\it
  isogenous} if there exists an isogeny $X \arrow Y$. In fact given an
isogeny $f: X \arrow Y$ there exists an isogeny $g: Y \arrow X$ such
that $gf=n_X \in \End(X)$ for some natural number $n$; in particular,
the relation `is isogenous to' is indeed an equivalence relation.

For an abelian variety $X$ we write $\End^0(X)$ to denote the tensor
product $\End(X) \otimes \Q$. Note that if $X$ and $Y$ are isogenous,
then $\End^0(X) \iso \End^0(Y)$ via pullback by the isogenies in
either direction. Therefore as a corollary of Theorem
\ref{thm-complete-reducubility} we get the following \cite[\S 19,
  Corollary 2]{Mumford1970}.

\begin{corollary} \label{corollary-completered}
Let $X$ be an abelian variety. If $X$ is simple, then $\End^0(X)$ is a
division algebra over $\Q$. For any abelian variety, if $X$ is
isogenous to $X_1^{n_1} \times \cdots \times X_k^{n_k}$ a product of
mutually non-isogenous simple abelian varieties, then
\begin{align*}
\End^0(X) = M_{n_1}(D_1) \times \cdots \times M_{n_k}(D_k)
\end{align*}
where $D_i$ is the division algebra $\End^0(X_i)$ and
$M_{n}(D_i)$ is the ring of $n \times n$ matrices over $D_i$.
\end{corollary}
So reduction to the case of simple abelian varieties is
straightforward. It is a much deeper problem to determine the
$\Q$-division algebras which can appear as $\End^0(X)$ for $X$ a
simple abelian varieties. The key point is that the endomorphism
algebra has finite rank and is equipped with a positive involution
\cite[\S 19 Corollary 3, \S 21 Theorem 1]{Mumford1970}:
\begin{proposition} \label{prop-rosati}
Let $X$ be a simple abelian variety. Then $D=\End^0(X)$ is a
finite-rank $\Q$-division algebra. Moreover $D$ carries an
involution $x \mapsto x'$ called the Rosati involution. This
involution is positive-definite in the sense that if $x \in D$ is any
nonzero element, then $\Tr_\Q(xx') > 0$, where $\Tr_\Q$ is the reduced
trace over $\Q$ of the division algebra $D$.
\end{proposition}
Now there is a classification (due to Albert) of all finite-rank
$\Q$-division algebras with a positive involution, and together with
some extra geometric restrictions this gives a complete list of
possibilities for $\End^0(X)$ when $X$ is a simple abelian
variety. Chapter 4 of Mumford's book gives an exposition of the
classification; here we only state a weak form of the result.

\begin{theorem} \label{thm-albertclass}
Let $X$ be a simple abelian variety and $D=\End^0(X)$ its endomorphism
algebra. Then $D \otimes_\Q \R$ is isomorphic as an $\R$-algebra with
involution to one of the following algebras:
\begin{itemize}
\item $\R \times \cdots \times \R$, with $x \mapsto x'$ the trivial
  involution;
\item $\H \times \cdots \times \H$ where $\H$ is the algebra of
  quaternions, and $x' = \overline{x}$, the usual conjugate;
\item $M_2(\R) \times \cdots \times M_2(\R)$ where $M_2(\R)$ is the
  algebra of $2 \times 2$ real matrices, and $x' = x^t$, the
  transpose;
\item $M_2(\C) \times \cdots \times M_2(\C)$, where $M_2(\C)$ is the
  algebra of $2 \times 2$ complex matrices, and $x' = x^*$, the
  conjugate transpose.
\end{itemize}
\end{theorem}
Combining this result with Corollary \ref{corollary-completered} one
can deduce the following description of $\End^0(X) \otimes \R$ for an
arbitrary abelian variety $X$.
\begin{corollary} \label{corollary-endoalgebra}
Let $X$ be an abelian variety. Then $\End^0(X) \otimes \R$ is
isomorphic as an algebra with involution to a product
\begin{align*}
\prod_i M_{r_i}(\R) \times \prod_j M_{s_j}(\C) \times \prod_k M_{t_k}(\H)  
\end{align*}
with involution given by conjugate transpose on each factor. The
bilinear pairing $\left<x, y \right>= \operatorname{Tr}(xy^*)$ defines
an inner product on $\End^0(X) \otimes \R$.
\end{corollary} 

Finally we must explain how the automorphism group $\Aut(X)$ sits
inside the algebra $\End^0(X) \otimes \R$. If $X$ is an abelian
variety, $\Aut(X)$ is the group of units $\End(X)^\times$ in the
endomorphism ring $\End(X)$. Furthermore $\End(X)$ is a lattice in the
vector space $\End^0(X) \otimes \R$, therefore induces a
$\Q$-structure on $\End^0(X) \otimes \R$, and this $\Q$-structure
determines $\End(X)$ as a subring of the $\R$-algebra $\End^0(X)
\otimes \R$ up to finite index. So from the previous corollary we have
the following:

\begin{corollary} \label{corollary-arithmetic}
Let $X$ be an abelian variety. Then $(\End^0(X) \otimes \R)^\times$ is
an algebraic group defined over $\Q$, and $\Aut(X)$ is an arithmetic
subgroup.
\end{corollary}

\section{The N\'eron--Severi space of an abelian variety}

In this final section, we explain how the N\'eron--Severi space of an
abelian variety can be identified with a subspace of the space
$\End^0(X) \otimes \R$. This identification allows us to describe the
action of the automorphism group on the ample cone in terms of
matrices, and applying the results of Section \ref{sect-reductionthy}
we get a proof of Theorem \ref{theorem-maintheorem}.

First we define the {\it N\'eron--Severi space} of an abelian variety
$X$ to be the finite-dimensional real vector space $N^1(X) :=
(\Pic(X)/\Pic^0(X)) \otimes \R$. Our first task is to identify
$N^1(X)$ with a subspace of $\End^0{X} \otimes \R$. By linearity, to
make such an identification it suffices to identify a Cartier divisor
with an element of $\End^0(X)$, which we do in the following way:
\begin{align*}
\xymatrix{ \Pic(X) \ar[r]^-{\phi} & \Hom(X,\widehat{X}) \otimes \Q
  \ar[r]^-{\psi} & \End^0(X) \\ D \ar@{|->}[r] & \phi_D \ar@{|->}[r] &
  \phi_L^{-1} \phi_D }
\end{align*}
The notation of the diagram is as follows. The variety $\widehat{X}$
is the {\it dual abelian variety} of $X$, which can be identified with
$\Pic^0(X)$. The homomorphism $\phi_D: X \arrow \widehat{X}$ is
defined by $\phi_D(x) = T_x^*(D) \otimes D^{-1}$, where $T_x$ is
translation by the point $x \in X$. Finally, $L$ is any (fixed) ample
line bundle on $X$; ampleness implies that $\phi_L$ is an isogeny $X
\arrow \widehat{X}$, and therefore has an inverse $ \phi_L^{-1} \in
\Hom(\widehat{X},X) \otimes \Q$. One checks that the kernel of $\phi$
is exactly the subgroup $\Pic^0(X)$ of numerically trivial line
bundles on $X$, and that $\psi$ is an isomorphism. Therefore tensoring
with $\R$ gives the claimed embedding $N^1(X)
\hookrightarrow \End^0{X} \otimes \R$. (See \cite{Mumford1970} for
proofs of all the assertions here.)

The automorphism group $\Aut(X)$ acts naturally on $N^1(X)$ via
pullback of divisors: an automorphism $f$ maps a divisor $D$ to
$f^*D$. Using the diagram above, we can extend this to an action of
$\Aut(X)$ on the whole algebra $\End^0(X) \otimes \R$. To compute the
action we need the following formula \cite[\S 15, proof of Theorem
  1]{Mumford1970}:
\begin{lemma} \label{lemma-matrixaction}
Let $f: X \arrow Y$ be an isogeny of abelian varieties, with dual
isogeny $\widehat{f}: \widehat{Y} \arrow \widehat{X}$. Let $D$ be a
line bundle on $Y$. Then $\phi_{f^*D} = \widehat{f} \circ \phi_D \circ
f$.\end{lemma} Via the diagram we can then work out the extension to
an action of $\Aut(X)$ on $\End^0(X) \otimes \R$: one computes that $f
\in \Aut(X)$ acts by the formula
\begin{align*}
f \cdot x = \phi_L^{-1} \circ \widehat{f} \circ
\phi_L \circ x \circ f
\end{align*}
 for any $x \in \End^0(X) \otimes \R$. Moreover the same formula
 defines an action of the whole group of units $(\End^0(X) \otimes
 \R)^\times$ on $\End(X)^0 \otimes \R$, and this action fixes the
 subspace $N^1(X)$.
 
To complete the picture, we observe \cite[p.189]{Mumford1970} that the
map $x \mapsto \phi_L^{-1} \circ \widehat{x} \circ \phi_L$ is by
definition exactly the Rosati involution on $\End^0(X) \otimes \R$
mentioned in Proposition \ref{prop-rosati}. (The involution therefore
depends on the choice of an ample line bundle on $X$, but the
dependence does not matter for our purposes.) For an element $e
\in \End^0(X) \otimes \R$ let us denote $\phi_L^{-1} \circ \widehat{e}
\circ \phi_L$ by $e'$, so that the action of $(\End^0(X) \otimes
\R)^\times$ on $\End(X)^0 \otimes \R$ is now given by the formula $x
\mapsto f' \circ x \circ f$. By Theorem \ref{thm-albertclass} there is
an isomorphism of $\End(X)^0 \otimes \R$ with a product of matrix
algebras which takes the Rosati involution $e \mapsto e'$ to the
conjugate transpose involution $x \mapsto x^*$. Using this isomorphism
to translate the action above into matrix terms, we get the following
theorem.

\begin{theorem}
Let $X$ be an abelian variety. Then the $\Q$-algebraic group
$(\End^0(X) \otimes \R)^\times$ acts on $N^1(X)$ by the formula $F: D
\arrow F^*DF$, and this extends the action of $\Aut(X)$ on $N^1(X)$ by
pullback of line bundles.
\end{theorem}

At this point we have an explicit description in terms of matrices of
the action of the automorphism group on the N\'eron--Severi space. To
complete the proof of the cone conjecture, we need to identify the
ample cone in the same terms (i.e. as a cone in a space of matrices),
and then apply the results of reduction theory from Section
\ref{sect-reductionthy}.

The first step is to identify the image of our embedding $N^1(X)
\hookrightarrow \End^0(X) \otimes \R$, viewing the target as a product
of matrix algebras as in Corollary \ref{corollary-endoalgebra}. Again
the key here is the Rosati involution: as we have just seen, in matrix
terms this is simply conjugate-transposition $x \mapsto x^*$.

\begin{theorem} \label{theorem-amplecone}
Let $X$ be an abelian variety. Then $N^1(X) \subset \End^0(X) \otimes
\R$ is exactly the fixed subspace of the Rosati involution. If
$\End^0(X) \otimes \R$ is isomorphic to a product of matrix algebras
\begin{align*}
\prod_i M_{r_i}(\R) \times \prod_j M_{s_j}(\C) \times \prod_k
M_{t_k}(\H)
\end{align*}
then $N^1(X)$ is isomorphic to the subspace 
\begin{align*}
\bigoplus_i \mathcal{H}_{r_i}(\R) \oplus \bigoplus_j
\mathcal{H}_{s_j}(\C) \oplus \bigoplus_k \mathcal{H}_k(\H)
\end{align*}
where $\mathcal{H}_r$ denotes the space of $r \times r$ symmetric or
Hermitian matrices. Moroever, the ample cone $\Amp(X)$ is the direct
sum of the positive-definite cones $P_r(k)$ in each of the direct
summands $\mathcal{H}_r(k)$ of $N^1(X)$.
\end{theorem} 
We use additive notation for $N^1(X)$ to emphasise the point that it
need not be a subalgebra of $\End^0(X) \otimes \R$: for divisors $D_1$
and $D_2$ fixed by the Rosati involution, their product
$D_1D_2~\in~\End^0(X) \otimes \R$ need not be fixed, as one can check
using suitable matrices.

Theorem \ref{theorem-amplecone} shows in particular that the ample
cone of an abelian variety is a self-dual homogeneous cone, since it
is a direct sum of cones on the list of Theorem
\ref{theorem-homogcones}. Indeed since it is a sum of cones of the
form $P_r(k)$, Theorem \ref{theorem-vinberg-auts} tells us that
$G(\Amp(X))^0$ acts transitively on $\Amp(X)$. Moreover from Theorem
\ref{theorem-amplecone} together with the description of the
automorphism group in Theorem \ref{theorem-vinberg-auts} we can also
deduce the following:
\begin{corollary} \label{theorem-surjection}
Let $X$ be an abelian variety. Then the homomorphism 
\begin{align*}
\left( \End^0(X) \otimes \R \right)^\times &\arrow G(Amp(X))^0 \\
M &\mapsto \left( D \mapsto M^*DM \right)
\end{align*}
given by the action of $(\End^0(X) \otimes \R)^\times$ on $N^1(X)$ is
surjective.
\end{corollary}
{\bf Proof:} Suppose for simplicity that $N^1(X)$ has a single direct
summand, say $N^1(X)=\mathcal{H}_{r}(\R)$. By Theorem
\ref{theorem-vinberg-auts} the identity component $G(\Amp(X))^0$ of
the automorphism group of the ample cone is exactly the group of
linear transformations of $N^1(X)$ of the form $D \mapsto M^*DM$ with
$M \in GL(r,\R)$. Such a linear transformation is the image under the
homomorphism of $M \in GL(r,\R) = (\End^0(X) \otimes \R)^\times$,
which proves surjectivity. The proof in the case of more than one
direct summand works in the same way, since by Theorem
\ref{theorem-homogcones} the identity component $G(\Amp(X))^0$ is the
direct product of the identity components of the automorphism groups
of the direct summand cones. QED Corollary \ref{theorem-surjection}

We can now complete the proof of Theorem \ref{theorem-maintheorem}.
The inclusion $N^1(X) \subset \End^0(X) \otimes \R$ endows $N^1(X)$
with a $\Q$-structure, and by Corollary \ref{theorem-surjection} the
connected component $\mathcal{G}(\Amp(X))^0$ is a $\Q$-algebraic
subgroup of $GL(N^1(X))$. By Corollary \ref{corollary-arithmetic} the
automorphism group $\Aut(X)$ is an arithmetic subgroup of $(\End^0(X)
\otimes \R)^\times$, so the image $\Gamma$ of $\Aut(X)$ in
$GL(N^1(X))$ is an arithmetic subgroup of $\mathcal{G}(\Amp(X))^0$.
Therefore by Theorem \ref{theorem-ash} and \ref{theorem-looijenga}
there is a rational polyhedral fundamental domain for the action of
$\Gamma$ on $\Amp(X)_+=\Nef{X}^e$ as required. QED Theorem
\ref{theorem-maintheorem}

\end{document}